\documentclass{amsart}

\usepackage{times,epsf}

\begin{document}

\newtheorem{thm}{Theorem}[section]
\newtheorem{lem}[thm]{Lemma}
\newtheorem{cor}[thm]{Corollary}

\theoremstyle{definition}
\newtheorem{defn}{Definition}[section]

\theoremstyle{remark}
\newtheorem{rmk}{Remark}[section]

\def\square{\hfill${\vcenter{\vbox{\hrule height.4pt \hbox{\vrule
width.4pt height7pt \kern7pt \vrule width.4pt} \hrule height.4pt}}}$}
\def\T{\mathcal T}

\newenvironment{pf}{{\it Proof:}\quad}{\square \vskip 12pt}

\title{Minimal Planes in Hyperbolic Space}
\author{Baris Coskunuzer}
\address{Department of Mathematics \\ Princeton University \\ Princeton, NJ 08544}
\email{baris@math.princeton.edu}

\maketitle


\newcommand{\cirD}{\overset{\circ}{D}}

\newcommand{\Si}{S^2_{\infty }}

\newcommand{\BH}{\Bbb H}
\newcommand{\BR}{\Bbb R}
\newcommand{\BC}{\Bbb C}
\newcommand{\BZ}{\Bbb Z}

\begin{abstract}

In this paper we show a generic finiteness result for least area
planes in $\BH^3$. Moreover, we prove that the space of minimal
immersions of disk into $\BH^3$ is a submanifold of product bundle
over a space of immersions of circle into $\Si(\BH^3)$ and the
bundle projection map is when restricted to this submanifold is
Fredholm of index zero. Using this, we also show that the space of
minimal planes with smooth boundary curve at infinity is a
manifold.

\end{abstract}

\section{Introduction}

The asymptotic Plateau problem in $\BH^3$ has been solved by
Anderson in [A1], by using geometric measure theory methods. In
this paper, we will consider this problem, by using geometric
analysis, and global analysis methods. Indeed, in a sense, we will
try to translate Anderson's results to these two fields to get
another perspective. We do not have any extra results about the
asymptotic Plateau problem, but we got a very nice global picture
of space of minimal planes, and by using this, we will prove
generic finiteness results.

A minimal plane P is a plane such that the mean curvature H is 0
at every point of x, i.e. $\forall x \in P$, $H(x)=0$. By Eells
and Sampson's results in [ES], for any minimal plane $\Sigma$ in
$\BH^3$, there is a conformal harmonic parametrization map $u:D^2
\rightarrow \BH^3$ with $u(D^2)=\Sigma$. This establishes the
starting point for our approach.

In this paper, we will only consider minimal planes with
$C^{3,\mu}$ regular asymptotic boundary. We will consider this
space as the space of minimal immersions of $D^2$ into $\BH^3$,
i.e. conformal harmonic maps from disk to hyperbolic space. Then,
one can think this space as a subspace of harmonic maps from disk
to hyperbolic space, and the minimal immersions in this space are
the conformal ones. As a result, we study the minimal maps as a
subspace of harmonic maps.

The space of harmonic maps from disk to hyperbolic space has very
nice features by Li and Tam's results. With some regularity
conditions, the space of harmonic maps from $D^2$ to $\BH^3$ can
be identified with their induced boundary parametrization
$\partial D^2$ to $\Si(\BH^3)$, by existence and uniqueness of
harmonic extensions in [LT1] and [LT2]. So, the space of minimal
maps can be considered as the subspace of immersions of $S^1$ into
$\Si(\BH^3)$.

Next, we will show that the space of minimal maps is a submanifold
of the space of immersions of $S^1$ into $\Si(\BH^3)$. Also, by
factoring out the parametrizations, we prove that the space of
minimal planes is a manifold, too.

Moreover, we will prove that the boundary restriction map from the
space of minimal immersions into the immersions of $S^1$ into
$\Si(\BH^3)$ is Fredholm map of index 0. In other words, for a
generic $ C^{3,\mu}$ Jordan curve $\Gamma \subset \Si(\BH^3)$, the
space of minimal planes spanning this curve, $M(\Gamma)$, is a
collection of isolated points.

\textbf{Theorem A:} $M=\{(\alpha,u)\in C^{3,\mu}(S^1,S^2)\times
C^3(S^1,S^1) |\  \widetilde{\alpha \circ u}(D^2) \mbox{ is
minimal} \}$ is a submanifold of the product bundle
$C^{3,\mu}(S^1,S^2)\times C^3(S^1,S^1)$, and the bundle projection
map when restricted to $M$, $\Pi_1|_M$, is Fredholm of index 0.

\textbf{Theorem B:} $M_0=\{\Sigma\subset \BH^3 | \Sigma \mbox{
minimal, }\partial_\infty\Sigma\subset\Si(\BH^3) \mbox{ is class
of }C^3\}$ the space of immersed minimal planes is a manifold.

The other main result is the generic finiteness of least area
planes spanning a fixed curve $\Gamma \subset \Si(\BH^3)$. By
Theorem A, for a generic curve $\Gamma$, the set of least area
planes spanning it, $M(\Gamma)$, is a collection of isolated
points. So all we need to show is that $M(\Gamma)$ is compact. We
will accomplish this by showing that this set is a subset of a
compact set by using boundary regularity results on minimal planes
in hyperbolic space by Hardt and Lin.

\textbf{Theorem C:} Let $\Gamma \subset \Si(\BH^3)$ be a
$C^{3,\mu}$ generic Jordan curve. Then there are finitely many
least area planes spanning this curve.

The organization of the paper is as follows: In Section 2, we will
give some preliminary definitions and theorems, which will be used
for the rest of the paper. In Section 3, we study the structure of
the space of minimal planes, and show that it is a manifold. In
Section 4, we will prove the generic finiteness result. Finally,
in Section 5, we will have some concluding remarks.

\subsection{Acknowledgements:}

I am very grateful to my advisor David Gabai for his continuous
encouragement and advice. I would like to thank Simon Brendle,
Alice Chang, Paul Yang, Fengbo Hang and Fang-Hua Lin for very
helpful conversations.

\section{Preliminaries}

In this section, we will parametrize the space of minimal
immersions of disk into hyperbolic space.

\begin{lem}$[$ES$]$(Minimal immersions and harmonic maps)

if $\varphi$ is an isometric immersion, then $\varphi$ is minimal
if and only if it is harmonic.

\end{lem}

\begin{thm}

Let $\Sigma$ be a plane in $\BH^3$. Then $\Sigma$ is minimal if
and only if there is a conformal harmonic map $u:D^2 \rightarrow
\BH^3$ with $u(D^2)=\Sigma$.

\end{thm}

\begin{pf}
Let $\Sigma$ be a minimal plane in $\BH^3$. Let $u:D^2 \rightarrow
\BH^3$ be an immersion with $u(D^2)=\Sigma$. Let $g$ be the
induced metric on $D^2$. Then $u:(D^2,g) \rightarrow \Sigma
\subset \BH^3$ is an isometric immersion. By Lemma 2.1, we
conclude that u is harmonic with respect to the induced metric
$g$.

By uniformization theorem, $(D^2,g)$ is conformally equivalent to
$\BC$ or Euclidean disk $(D^2,g_0)$. But, since g is induced by
$\Sigma$, mean curvature 0 plane, in $\BH^3$, constant sectional
curvature -1 space, then $(D^2,g)$ is conformally equivalent to
$(D^2,g_0)$. So, there exist a conformal map $\phi
:(D^2,g_0)\rightarrow (D^2,g)$. Since harmonicity respects
conformality, $u \circ\phi : (D^2,g_0) \rightarrow \Sigma$ is
harmonic. Then, there exist a conformal harmonic map $v: D^2
\rightarrow \Sigma$.

The reverse direction is trivial, as by definition, for any
conformal harmonic map, the mean curvature is 0 at every point in
the image.
\end{pf}

So, by using this theorem, we can consider minimal planes as
conformal harmonic maps. Now, let's focus on harmonic maps.

\begin{lem}$[\mbox{LT1}]$$[\mbox{LT2}]$(Harmonic extensions)

If $\gamma:S^1_\infty(\BH^2) \rightarrow \Si(\BH^3)$ a $C^1$
immersion, then there exist unique $C^1$-regular harmonic
extension map $\overline{\BH^2} \rightarrow \overline{\BH^3}$.

\end{lem}

\begin{lem}$[$L$]$$[$To$]$(Boundary regularity of minimal planes)

If $\Gamma \subset \Si(\BH^3)$ is a $C^1$ curve, and $\Sigma$ is
minimal plane with $\partial_\infty\Sigma = \Gamma$, then $\Sigma
\cup \Gamma$ is $C^1$, too.

\end{lem}

\begin{thm}

Let $M=\{u:D^2\rightarrow \BH^3 | u$ minimal immersion and $
u|_{\partial D^2}$ is $C^1 \}$, and let $A=\{\alpha:S^1
\rightarrow S^2 | \alpha$ immersion $\}$. Then $\beta: M
\rightarrow A$ boundary restriction map is one to one.

\end{thm}

\begin{pf}
Let u be a minimal immersion with $ u|_{\partial D^2}$ is $C^1$.
Then $\alpha := u|_{\partial D^2}$ is in A, the image of u under
the boundary restriction map. This map is 1-1, since by Lemma 2.4,
the image of u is $C^1$ and by theorem 2.3, there exist unique
$C^1$ harmonic extension of $\alpha$, which must be u. So, $\beta$
embeds M into A.
\end{pf}

So, we can identify M with $\beta(M) \subset A$, and we can
consider the space of minimal immersions M as a subspace of space
of immersions of $S^1$ into $S^2$. In the following sections, we
will use this identification.

\section{The Space of Minimal Planes}

In this section, we will show that the space of minimal immersion
is a manifold, and boundary restriction map is Fredholm of index
0. Basically, we will adapt the techniques of Tomi and Tromba in
[TT] and [To].

We will start by defining the following spaces: Let $\mu \geq 0$.

$A= \{\alpha \in C^{3,\mu}(S^1,S^2)| \alpha$ immersion $\}$

$D=\{u \in C^3(S^1,S^1)| u$ diffeomorphism and satisfies three
point condition, i.e. $u(e^{\frac{2}{3} k \pi i}) = e^{\frac{2}{3}
k \pi i}, k=1,2,3 \}$

$M=\{f:D^2\rightarrow \BH^3 | f(D^2)$ minimal and $ f|_{\partial
D^2} \in A \}$

$N=A \times D$

Here, A is open subset of $C^{3,\mu}(S^1,S^2)$, D is a smooth
manifold, and so is N.

Define the conformality operator $k:N \rightarrow C^2(S^1,\BR)$
such that $$k(\alpha,u)= \frac{\partial}{\partial r}(\alpha \circ
u) \cdot \frac{\partial}{\partial \theta}(\alpha \circ u)$$

Here, $r,\theta$ represents polar coordinates on $D^2$ and
functions $\alpha \circ u$ identified with their harmonic
extensions, $\widetilde{\alpha \circ u}$. So, conformality
operator is dot product of partial derivatives of $\alpha \circ
u$.

\begin{lem}

$ker(k) \subset N$ consists of conformal harmonic maps, i.e.
$k(\alpha,u)=0$ if and only if the harmonic extension
$\widetilde{\alpha \circ u}$ is conformal.

\end{lem}

\begin{pf}
Let $$\delta= (\cos\theta - i \sin\theta)(\frac{\partial}{\partial
r} - \frac{i}{r}\frac{\partial}{\partial \theta})$$ be the complex
differential. Then $$k(\alpha,u)= - \mbox{Im}(z^2
\delta(\alpha\circ u )\cdot \delta(\alpha\circ u ))$$ But
$\delta(\widetilde{\alpha \circ u} )\cdot \delta(\widetilde{\alpha
\circ u})$ is the Hopf differential of $\widetilde{\alpha \circ
u}$, $$Q(\widetilde{\alpha \circ u})= (|\widetilde{\alpha \circ
u}|^2_x - |\widetilde{\alpha \circ u}|^2_y) - 2i
((\widetilde{\alpha \circ u})_x \cdot (\widetilde{\alpha \circ
u})_y)$$ The Hopf differential, $Q$, is holomorphic for harmonic
maps, [ER]. Then, $Q(\widetilde{\alpha \circ u}) = 0$ if and only
if $\widetilde{\alpha \circ u}$ is conformal. Since Im$(z^2
\delta(\widetilde{\alpha \circ u} )\cdot \delta(\widetilde{\alpha
\circ u} ))$ is harmonic function, then $k(\alpha,u)= 0$ implies
$Q(\widetilde{\alpha \circ u}) = 0$, which means
$\widetilde{\alpha \circ u}$ is conformal.
\end{pf}

To understand the structure of ker($k$), we will study the
derivative of $k$:

$$D_\alpha k(\alpha,u)\langle \beta\rangle = \frac{\partial}{\partial r}(\beta \circ
u) \cdot \frac{\partial}{\partial \theta}(\alpha \circ u) +
\frac{\partial}{\partial r}(\alpha \circ u) \cdot
\frac{\partial}{\partial \theta}(\beta \circ u)$$

$$D_u k(\alpha,u)\langle v\rangle= \frac{\partial}{\partial r}\left[(D\alpha)\circ u
\langle v\rangle\right] \cdot \frac{\partial}{\partial
\theta}(\alpha \circ u) + \frac{\partial}{\partial r}(\alpha \circ
u)\cdot \frac{\partial}{\partial \theta}\left[(D\alpha)\circ u
\langle v\rangle\right]$$

\begin{lem}$[$ADN$]$ Let $a_1,a_2,a_3,b \in C^2(S^1,\BR^3)$ such that
for every $p\in S^1$, the vectors $a_1(p),a_2(p)$, and $a_3(p)$
are linearly independent. Then the following system is an elliptic
boundary value system.

$\ \ \ \ \ \ \ \ \ \ \ \ \ \ \ \ \ \ \ \ \ \ \triangle w \ \ \ \ \
=\ \ \ \  f \ \ \ \ \ \ \ \  \mbox{on} \ \ D^2 $

$\left.
\begin{array}{ccc}
\ \ \ \ \ \ \ \ \ \ \ \ \ \ \ \ \ \    a_1 \cdot w & = & g_1 \\
\ \ \ \ \ \ \ \ \ \ \ \ \ \ \ \ \ \    a_2 \cdot w & = & g_2 \\
 a_3 \cdot \frac{\partial w}{\partial r} + b \cdot
\frac{\partial w}{\partial \theta} & = & g_3  \\
\end{array}
 \right\rbrace \ \  \mbox{on}\ \  S^1$

\end{lem}

\begin{lem}$[$BT$]$
Let $F, X_1, X_2$ be Banach spaces and let $A_k: F \rightarrow
X_k$ with $k=1,2$ be linear operators. Let $A:=(A_1,A_2): F
\rightarrow X_1 \times X_2$ be Fredholm operator. Then
$L:=A_1|_{A_2^{-1}(0)}$ is also a Fredholm operator.

\end{lem}

\begin{cor}
The map $ w \rightarrow a_3 \cdot \frac{\partial w}{\partial r} +
b \cdot \frac{\partial w}{\partial \theta}$ from $\{w\in
C^{3,\mu}(S^1,\BR^3) | \triangle_{\BH} w=0, a_1 \cdot w= a_2 \cdot
w= 0\}$ into $C^2(S^1, \BR)$ is a Fredholm operator.
\end{cor}

\begin{pf}
Since the system in Lemma 3.2 is elliptic boundary value system,
then the operator
$$w \rightarrow (\triangle w, a_1 \cdot w,a_2
\cdot w,a_3 \cdot \frac{\partial w}{\partial r} + b \cdot
\frac{\partial w}{\partial \theta})$$

for the suitable spaces is Fredholm, [LRW]. Then by Lemma 3.3, the
map
$$ w \rightarrow a_3
\cdot \frac{\partial w}{\partial r} + b \cdot \frac{\partial
w}{\partial \theta}$$ from $ \{w\in C^3(S^1,\BR^3) | \triangle
w=0, a_1 \cdot w= a_2 \cdot w= 0\} \mbox{ into } C^2(S^1, \BR)$ is
also Fredholm. Since $\triangle w=0$ condition is essentially same
in the interior with $\triangle_{\BH} w=0$ condition, then the
result follows.

\end{pf}

\begin{thm}
$D_u k(\alpha,u): T_u D \rightarrow C^2(S^1,\BR)$ is Fredholm
where $(\alpha,u)\in N$.
\end{thm}

\begin{pf}
For $v\in T_u D$,
$$D_u k(\alpha,u)\langle v\rangle= \frac{\partial}{\partial r}[(D\alpha)\circ u
\langle v\rangle] \cdot \frac{\partial}{\partial \theta}(\alpha
\circ u) + \frac{\partial}{\partial r}(\alpha \circ u)\cdot
\frac{\partial}{\partial \theta}[(D\alpha)\circ u \langle
v\rangle]$$

Let $N^\alpha=\{\alpha \circ u | u\in D\}$. Then by using the
isomorphism $T_u D \rightarrow T_{\alpha \circ u}N^\alpha$ with $v
\rightarrow (D\alpha)\circ u \langle v\rangle$, we can simplify
$D_u k(\alpha,u)$ with $w=(D\alpha)\circ u \langle v\rangle$. By
abuse of notation, we will have $D_u k(\alpha,u): T_{\alpha \circ
u}N^\alpha\rightarrow C^2(S^1,\BR)$

$$D_u k(\alpha,u)\langle w\rangle = a_3 \cdot \frac{\partial w}{\partial r} + b
\cdot \frac{\partial w}{\partial \theta}$$ where
$a_3=\frac{\partial}{\partial \theta} (\alpha \circ u)$ and
$b=\frac{\partial}{\partial r} (\alpha \circ u)$.

Then, we can find vector functions $a_1,a_2 \in C^2(S^1,\BR^3)$
such that $\{a_1,a_2,a_3\}$ form an orthogonal triple and
$$T_{\alpha \circ u} N^\alpha = \{w\in
C^3(S^1,\BR^3) | \triangle_{\BH} w=0,\  a_1 \cdot w= a_2 \cdot w=
0\}$$ by the definition of $a_3$. By Lemma 3.2 and Lemma 3.3, $D_u
k$ is a Fredholm operator.

\end{pf}

\

\begin{thm}
The Fredholm index of $D_u k(\alpha,u)$ is 0 for $(\alpha,u)\in
N$.

\end{thm}

{\em Proof:} We will prove this theorem in 3 steps.

\textbf{Claim 1:} Let $$Z=\{h \in C^2(S^1,\BR) | \int_{S^1}h
d\theta = \int_{S^1} h \cos\theta d\theta = \int_{S^1}h\sin\theta
d\theta= 0\}$$ be the subspace of $C^2(S^1,\BR)$. Then the image
of $k$ is contained in Z, i.e. $k(A\times D) \subset Z \subset
C^2(S^1,\BR) $.

\begin{pf}
By above, $k(\alpha \circ u)= - Im(z^2 Q(\widetilde{\alpha \circ
u})|_{\partial D^2})$ Since $Q$ is holomorphic, the result follows
by Cauchy's theorem.
\end{pf}

\textbf{Claim 2:} $D_u k(id,id): T_{id} D \rightarrow Z$ is an
isomorphism, where $id$ is the identity of $S^1$.

\begin{pf}
First, we will show that $D_u k$ is injective. If we identify the
elements of $D$ with their harmonic extensions, then the elements
of $ T_{id} D$ will be harmonic functions from $D^2$ to $\BR^2$.
so, we can think of them as complex valued functions. Let $v\in
T_{id} D$.

\begin{eqnarray*}
  D_u k(id,id)\langle v\rangle & = & \frac{\partial v}{\partial r} \cdot \frac{\partial id}{\partial \theta} + \frac{\partial id}{\partial r}\cdot \frac{\partial v}{\partial \theta}\\
   & = & Re(\frac{\partial v}{\partial r} \cdot \frac{\overline{\partial id}}{\partial \theta} + \frac{\partial v}{\partial \theta}\cdot \frac{\overline{\partial id}}{\partial r})\\
   & = & Re((z\delta v + \bar{z}\bar\delta v)(-i\bar{z}) +(iz\delta v - i\bar{z}\bar\delta v)\bar{z}) \\
   & = & Re(-i\bar{z}^2\bar\delta v)=Re(iz^2\delta\bar{v}).\\
\end{eqnarray*}

Since v is harmonic, $\delta \bar{v}$ is holomorphic. So, Re$(i
z^2 \delta\bar{v})=0$ implies $i z^2 \delta\bar{v}=0$ and v is
itself holomorphic. Now, since  for any $u\in D$, $u(S^1)=S^1$ and
$v\in T_{id} D$, v is tangential to $S^1$, then the function
$$w=\frac{v}{iz}$$ is real valued on $S^1$. By reflection principle,
we can extend w to a meromorphic function $\hat w$ on $\bar \BC$
with simple poles at $0$ and $\infty$. Then by elementary complex
analysis, $$w(z)=a + bz + \frac{\bar b}{z}$$ and $$v(z)=i(\bar b +
az + bz^2)$$ where $a \in \BR$ and $b\in \BC$. Then, by the three
point condition, v must be 0. So, $D_uk (id,id)$ is injective.

Now, we will show surjectivity. Let $h \in Z$. Then there exist a
holomorphic function $g$ with boundary values of class $C^2$ such
that $h=$Re$g$, $g(0)=g'(0)=0$. So, g can be written as $g=iz^2 f$
where $f$ is a holomorphic function with boundary values of class
$C^2$. We claim that there is a $v\in T_{id} D$, such that Re$(i
z^2 \delta\bar{v})=h$.

Set $$v= \bar{F} + H$$ where $F$ and $H$ are holomorphic and
$F'=f$. Obviously, for such $v$, Re$(i z^2 \delta\bar{v})=h$
holds. All we need to do is to find a suitable $H$ such that
$v=\bar{F} + H \in T_{id}D$. $v$ is tangential to $S^1$ if and
only if Re$(\bar z v)=0$ on $S^1$. So, we need

$$\begin{array}{ccc}
\mbox{Re}(\bar z H)= - \mbox{Re}(zF)& \mbox{or} & \mbox{Re}(\frac{H}{z})= - \mbox{Re}(zF)\\
\end{array}$$

Then, if we put
$$H= -z^2 F +i(\bar b + az +bz^2)$$ where $a\in
\BR$ and $b\in \BC$ are chosen in such a way that v satisfies the
three point condition. So, $D_u k(id,id)$ is surjective.
\end{pf}

\textbf{Claim 3:}For $(\alpha,u)\in N$, the Fredholm index of $D_u
k(\alpha, u )$ is 0.

\begin{pf}
Consider the map  $\Phi: N \rightarrow L(E_1,E_2)$ from $N$ into
bounded linear operators, with $\Phi(\alpha,u)=D_u k(\alpha,u)$
where $E_1=T_{id} D$ and $E_2=C^2(S^1,\BR)$. Here we identified
$T_u D$ with $T_{id} D$ for any $u\in D$, with abuse of notation.
Then $\Phi(N) \subset F(E_1,E_2) \subset L(E_1,E_2)$ where
$F(E_1,E_2)$ represents Fredholm operators. By [Sm], index is
continuous on Fredholm operators. Since for any $(\alpha,u)\in N$,
there exist a path in $N$ connecting $(\alpha,u)$ to $(id,id)$, by
the existence of homotopies in the corresponding spaces. Then,
continuity of index implies ind$(D_u k(\alpha,u))= $ ind$(D_u
k(id,id))=0$.
\end{pf}

Let $M=\{(\alpha,u) \in N | k(\alpha,u)=0,
\delta(\widetilde{\alpha\circ u})\neq 0\}$ represents the space of
minimal immersions of $D^2$ into $\BH^3$. Here, for any
$(\alpha,u) \in M$, $\widetilde{\alpha\circ u}: D^2 \rightarrow
\BH^3$ is minimal immersion.

\begin{thm}
M is a submanifold of $N$ with tangent space $T_{(\alpha,u)} M=
ker (Dk(\alpha,u))$.

\end{thm}

\begin{pf}
We will prove this theorem by showing that if $(\alpha,u)\in M
\subset N$ then $Dk(\alpha,u): T_{(\alpha,u)}N \rightarrow Z$ is
onto.

$$D k(\alpha,u)\langle \beta, v\rangle= \frac{\partial w}{\partial r}
\cdot \frac{\partial}{\partial \theta}(\alpha \circ u) +
\frac{\partial}{\partial r}(\alpha \circ u) \cdot \frac{\partial
w}{\partial \theta}$$

where $w=\beta\circ u + (D\alpha)\circ u \langle v\rangle$.

Consider the simple functional analytic fact: Let $X,Y$ be Banach
spaces, $X_0$ be a dense subspace of $X$, and $T:X\rightarrow Y$
be linear and continuous. If $T$ is onto and $T(X_0)$ is closed,
then $T(X_0)=Y$.

Set $X_0=$Range$\{\beta\rightarrow\beta\circ u\}\times T_u D$,
$X=C^3(S^1,TS^2)\times T_u D$, and $Y=Z$. We know
Range$(Dk(\alpha,u))$ is closed by Theorem 3.5. So, by the above
fact, it suffices to solve the equation $$\frac{\partial
w}{\partial r} \cdot \frac{\partial}{\partial \theta}(\alpha \circ
u) + \frac{\partial}{\partial r}(\alpha \circ u) \cdot
\frac{\partial w}{\partial \theta}=h$$ in $X$ for $h\in Z\subset
C^2(S^1,\BR)$.

In complex notation, this equation becomes Im$(z^2\delta
w\cdot\delta(\alpha\circ u))=$ Im$(z^2 g)$, where $g$ is
holomorphic in the unit disc and has boundary values $C^2(S^1)$.
Let's denote this space by $\Psi$.

Define another space $\Phi = \{ V:D^2 \rightarrow \BR^3 |
\nabla_{\bar z} V= 0,$ and $V|_{\partial{D^2}}$ is class of $C^2
\}$ where $\nabla$ is covariant derivative. Recall that $\varphi$
is harmonic map if and only if $\nabla_{\bar{z}}\partial_z
\varphi=0$.

Then it is sufficient to solve the equation $W\cdot F=g$, where
$g\in \Psi$ , $F=\delta(\widetilde{\alpha\circ u})$ are given, $F$
does not vanish in the unit disc and $W=\delta\widetilde{w}$ is
unknown. Since $\widetilde{\alpha\circ u}, \widetilde{w}$ are
harmonic maps, then
$$\nabla_{\bar{z}}\delta(\widetilde{\alpha\circ
u})=\nabla_{\bar{z}}\delta\widetilde{w}=0$$ so $F,W\in\Phi$.

Let $V_1,V_2\in\Phi$, then consider $$\partial_{\bar{z}}(V_1\cdot
V_2)=(\nabla_{\bar{z}}V_1\cdot V_2) +(V_1\cdot\nabla_{\bar{z}}V_2)
=0$$ so $V_1\cdot V_2 \in \Psi$.

We will show that there exist a solution to the equation $W\cdot
F=g$ by using the fact that Banach space $\Psi$ is a topological
algebra. Define the set $J= \Phi \cdot F$. If the equation is not
solvable for some $g\in\Psi$, then J is an ideal of $\Psi$. Since
any ideal is contained in some maximal ideal, we may conclude from
the Gelfand-Mazur Theorem that J is contained in the kernel of
some algebra homomorphism $\tau: \Psi \rightarrow \BC$.

By using the norm $$\|f\|=  {\left( \sum_{j=0}^{\infty}
(1+k^2){|a_k|}^2 \right)}^{\frac{1}{2}}$$

for $f=\sum a_k z^k \in\Psi$, we get
$$|\tau(id)|=\sqrt[n]{|\tau(id^n)}\leq\sqrt[n]{\|\tau\|}\sqrt[n]{\|id^n\|}$$

Then, $$|\tau(id)|\leq\limsup_{n\rightarrow
\infty}\sqrt[n]{\|id^n\|}=\limsup(1+n^2)^\frac{1}{2n}=1$$

So, $\xi_0=\tau(id)$ is a point of the unit disc and, for an
arbitrary $f=\sum a_k z^k \in\Phi$, we have $$\tau(f)=\sum
a_k\tau(id)^k=f(\xi_0)$$ Since $J\subset ker(\tau)$, for any
$W\in\Phi$, we have $$0=\tau(W\cdot F)= W(\xi_0)\cdot F(\xi_0)$$
But this implies $F(\xi_0)=0$, which contradicts to the immersion
assumption.
\end{pf}

\begin{lem}$[$TT$]$
 Let $X_1, X_2, Y$ be Banach spaces and
$L:X_1\times X_2 \rightarrow Y$ be a linear and surjective
operator such that $L_2:=L(0,\cdot):X_2\rightarrow Y$ is Fredholm.
Let $\Pi$ denote the projection $X_1\times X_2 \rightarrow X_1$
and $M:= ker L$. Then we have $ker(\Pi|_M)=ker(L_2)$, and
$coker(\Pi|_M)=coker(L_2)$, i.e. $\Pi|_M$ is Fredholm and
$\mbox{ind}(\Pi|_M)= \mbox{ind}(L_2)$.
\end{lem}

\begin{thm}
The restriction of the projection map $\Pi: A \times D \rightarrow
A$ to the submanifold $M$, $\Pi|_M$ is onto and  $C^2$ Fredholm of
index 0.
\end{thm}

\begin{pf}
$\Pi|_M$ is onto, since for any $\alpha \in A$, $\Gamma=
\alpha(S^1)\subset \Si(\BH^3)$ is a $C^{3,\mu}$ curve, and there
is a minimal (indeed, least area) plane $\Sigma \subset \BH^3$
spanning $\Gamma$ by Anderson's results [A1]. Then by Theorem 2.2,
there is a conformal harmonic parametrization of $\Sigma$, say
$\varphi$. Since there exist $u\in D$ such that
$\varphi|_{\partial{D^2}} = \alpha\circ u$, then $(\alpha,u)\in
M$. This implies $\Pi_M$ is onto.

Let $(\alpha,u)\in M$. Then by using the above lemma with,
$X_1=T_\alpha A$, $X_2=T_u D$, $Y=Z$, $L=Dk(\alpha,u)$, we get
ind$(\Pi|_M)=$ind$(D_u k(\alpha,u))$. Then by Theorem 3.6, the
result follows.
\end{pf}

So far, we proved that the space of minimal immersions, $M$, is a
submanifold of $N$. But, in this space, for a given minimal plane
$\Sigma \subset \BH^3$, there are different corresponding
parametrizations. For example, if $\Sigma=(\widetilde{\alpha \circ
u})(D^2)$ then $\Sigma$ is represented by $(\alpha,u)\in M$. But,
for any $v\in D$, $(\alpha\circ v,v^{-1}\circ u)\in M$ represents
same minimal plane $\Sigma$.

All results in this section up to now are valid for any $\mu \geq
0$. Now, we will show for $\mu=0$, the space of minimal planes is
a manifold. Let $M_0=\{\Sigma\subset \BH^3 | \Sigma$ minimal,
$\partial_\infty\Sigma\subset\Si(\BH^3)$ is class of $C^3\}$ be
the space of minimal planes. By above discussion to get $M_0$, we
have to factor out the group actions from $M$.

\begin{lem}

$M_0= M_1 / C$, where $M_1 = M/D$ and $C$ is the space of
conformal diffeomorphisms of the unit disk.

\end{lem}

\begin{pf}
Consider the $D$ action on $M$. For any $v\in M$,
$\Psi_v:M\rightarrow M$ such that $\Psi_v(\alpha,u)=(\alpha\circ
v, v^{-1}\circ u)$. Obviously, $D$ action does not change the
induced harmonic extension $\widetilde{\alpha \circ u}$ as the
composition of $(\alpha\circ v \circ v^{-1}\circ u)= \alpha \circ
u$. So, $M_1= \{ [( \alpha,u)] | [(\beta,v)]=[( \alpha,u)]$ if
$\beta\circ v=\alpha\circ u \}$. In a sense, we eliminated the
artificial augmentation in the base space N by $D$. In other
words, we reduced the parametrizations of minimal immersions to
space $A$, since $M_1 \simeq M\cap (A\times\{id\})$.

Now, consider $C= \{ c:D^2 \rightarrow D^2 | c$ conformal
$c|_{\partial D^2}$ is class of $C^3\}$. Now consider the action
of C on $M_1$. Define for any $c\in C$, $\Phi_c:M_1\rightarrow
M_1$ such that $\Phi_c([(\alpha,u)])=[(\alpha\circ c,u)]$, or by
using the equivalence $M_1 \simeq M\cap (A\times\{id\})$, then
$\Phi_c((\alpha,id))=(\alpha\circ c,id)$. So, $M_1 /
C=\{([\alpha],id) | \alpha\sim\beta$ if $\exists c\in C,$
$\beta=\alpha\circ c\}$.

Now, if we can show that for any minimal plane $\Sigma\in M_0$,
there exist unique $([\alpha],id)$ with
$\Sigma=\tilde{\alpha}(D^2)$, then we are done. Existence is true
by section 2. Now, assume $\exists \alpha, \beta\in A$ such that
$\Sigma=\tilde{\alpha}(D^2)=\tilde{\beta}(D^2)$. Then
${\tilde\beta}^{-1}\circ \tilde\alpha:D^2\rightarrow D^2$ is a
conformal diffeomorphism of the disk. This implies
$\beta^{-1}\circ\alpha\in C$, and $([\alpha],id)=([\beta],id)$, so
the uniqueness follows.
\end{pf}

\begin{lem}$[$MO$]$
Let a group G acts freely on a manifold X, i.e. $\forall g\in
G-\{id\}, \forall x\in X$, $g(x)\neq x$. Then the quotient space
$X/G$ is a manifold.
\end{lem}

\begin{thm}
The space of minimal planes, $M_0$, is a manifold.
\end{thm}

\begin{pf}
First, we will show that $M_1$ is a manifold, by showing the
action of D on M is free. Let $v \in D$, and
$\Psi_v((\alpha,u))=(\alpha,u)$. Then $(\alpha\circ v,v^{-1}\circ
u)=(\alpha,u)$. But, this implies $v^{-1}\circ u=u$ and so $v=id$.
The action is free, and by lemma 3.11, $M_1$ is a manifold.

Now, if we can show the action of $C$ on $M_1$ is free, then we
are done. Let $c\in C$, and $(\alpha,id)\in ((A \times \{id\})\cap
M) \simeq M_1$, with $\Phi_c((\alpha,id))=(\alpha,id)$. Then
$(\alpha\circ c,id)=(\alpha,id)$. But, this implies $\alpha\circ
c=\alpha$, and $c=id$. So, the action is free, and $M_0=M_1 /C$ is
a manifold.
\end{pf}

\section{Generic Finiteness}

In this section, our aim is to prove that for a generic
$C^{3,\mu}$ simple Jordan curve $\Gamma \subset \Si(\BH^3)$, there
exist finitely many least area planes in $\BH^3$ spanning
$\Gamma$. From now on, we fix a $\mu > 0$.

Let $M=\{(\alpha,u) \in N | \widetilde{\alpha\circ u}$ is minimal
immersion$\}$ represents again the space of minimal immersions in
$N$. We proved in previous section that $M$ is a submanifold of
the product bundle $N=A\times D$ with $D \hookrightarrow A\times D
\overset{\Pi}{\rightarrow} A$. Moreover, the projection map
$\Pi|_M \rightarrow A$ is Fredholm of index 0.

\begin{thm}(Sard-Smale)$[$Sm$]$
Let $f: X\rightarrow Y$ be a Fredholm map. Then the regular values
of $f$ are almost all of Y, i.e except a set of first category.
\end{thm}

\begin{cor}$[$Sm$]$
Let $f: X\rightarrow Y$ be a Fredholm map. Then for any regular
value $y\in Y$, $f^{-1}(y)$ is a submanifold of X whose dimension
is equal to index($f$) or empty.
\end{cor}

\begin{thm}
For almost all $\alpha \in A$, the set $(\Pi|_M)^{-1}(\alpha)$ is
a collection of isolated points in $M$.
\end{thm}

\begin{pf}
Since $\Pi|_M$ is onto by Theorem 3.9, $(\Pi|_M)^{-1}(\alpha)$ is
not empty for any $\alpha \in A$. Since $\Pi|_M$ is $C^2$ Fredholm
of index 0, by the corollary, $(\Pi|_M)^{-1}(\alpha)$ is 0
dimensional submanifold of M, for almost all $\alpha \in A$. The
result follows.
\end{pf}

Now, let $M^\alpha:=(\Pi|_M)^{-1}(\alpha)= M \cap
(\{\alpha\}\times D)$. So far we have shown that $M^\alpha$ is a
collection of isolated points. If we can show it is also finite,
then we are done.

\begin{lem}(Meeks-Yau)$[$MY$]$
Let $M^3$ be a compact Riemannian three-manifold whose boundary is
mean convex and let $\Gamma$ be a simple closed curve in $\partial
M$ which is null-homotopic in M; then $\Gamma$ is bounded by a
least area disk and any such least area disk is properly embedded.
\end{lem}

\begin{thm}
Let $\Gamma \subset \Si(\BH^3)$ be a $C^{3,\mu}$ regular, simple
Jordan curve. Then, any least area plane spanning $\Gamma$ is
properly embedded.
\end{thm}

\begin{pf}
By the boundary regularity results of Hardt and Lin in [HL], we
can find a sufficiently large $N>0$ such that in Poincare Ball
Model for $\BH^3$, $\partial B_N(0) \cap \Sigma$ is a simple
closed curve and $\Sigma - B_N(0)$ is properly embedded, indeed a
graph over an annulus. Now, since $\Sigma\cap B_N(0)$ is also
properly embedded by Lemma 4.4 as $B_N(0)\subset \BH^3$ is convex,
then result follows.
\end{pf}

\begin{lem}
Let $\Gamma \subset \Si(\BH^3)$ be a $C^{3,\mu}$ curve. Then for
any embedding $\varphi: D^2 \rightarrow \BH^3$ with
$\varphi(\partial{D^2})=\Gamma$ and $\varphi(D^2)$ least area
plane, the $C^{3,\mu}$ norm of such minimal immersions are
uniformly bounded by a constant depending on $\Gamma$, i.e.
$||\varphi||_{C^{3,\mu}} < C(\Gamma)$.
\end{lem}

\begin{pf}
This follows from the results of Hardt and Lin in the articles
[HL] and [L] about the boundary regularity of least area planes in
hyperbolic space. In these papers, the authors proved that if
$\Gamma \subset \Si(\BH^3)$ be a $C^{3,\mu}$ curve, let $p\in
\Gamma$, then there exist a neighborhood of point
$N(p)\subset\bar{\BH^3}$, such that for any minimal plane $\Sigma
\subset \BH^3$ spanning $\Gamma$, if $\phi$ parametrizes
$\Sigma\cap N(p)$, then $||\phi||_{C^{3,\mu}} < C(p,\Gamma)$. Now,
by using this local result, they conclude the boundary regularity
of least area planes. Here, since interior regularity is already
known, by using this local estimates and using simple
transformations, we can reach a global estimate, such that any
conformal harmonic parametrization of a least area plane spanning
$C^{3,\mu}$ curve has $C^{3,\mu}$ norm uniformly bounded by
constant depends only on curve and independent of the plane and
its parametrization.
\end{pf}

\begin{thm}
Let $\Gamma \subset \Si(\BH^3)$ be a $C^{3,\mu}$ be a generic
curve as described above. Then there are finitely many least area
planes spanning this curve.
\end{thm}

\begin{pf}
Take a parametrization of $\Gamma$ in $A$, say $\alpha_0$. Then by
previous sections, we know that all least area planes spanning
$\Gamma$ have parametrizations in the form of
$\widetilde{\alpha_0\circ u}$ with $(\alpha_0,u)\in M^{\alpha_0}$.
So if we can show that $M^{\alpha_0}$ is finite, then the result
follows. Now, since this set is bounded in $C^{3,\mu}$ topology by
a constant $C(\Gamma)$ by Lemma 4.6, it is compact in some
$C^{3,\delta}$ topology where $0<\delta<\mu$. Now, if we can also
show $M^{\alpha_0}$ is a collection of isolated points in
$C^{3,\delta}$ topology, then this will imply $M^{\alpha_0}$ is
finite.

Let $A_\mu=\{\alpha \in C^{3,\mu}(S^1,S^2)| \alpha$ immersion
$\}$, and let $M_\mu \subset A_\mu \times D$ be the submanifold
corresponding minimal immersions. Similarly, define
$A_\delta=\{\alpha \in C^{3,\delta}(S^1,S^2)| \alpha$ immersion
$\}$, and $M_\delta$. Now, by section 3, $\Pi_\mu:M_\mu
\rightarrow A_\mu$, and $\Pi_\delta:M_\delta \rightarrow A_\delta$
are Fredholm maps of index 0. Clearly, $M_\mu \subset M_\delta$
and $A_\mu \times D \subset A_\delta \times D$. Moreover,
$\Pi_\mu=\Pi_\delta|_{M_\mu}$. Now, we claim that if $\alpha \in
A_\mu \subset A_\delta$ is a regular value for $\Pi_\mu$, then
$\alpha$ is also a regular value for $\Pi_\delta$. This will imply
$M^{\alpha_0}$ is a collection of isolated points in
$C^{3,\delta}$ topology, by Lemma 4.3, and the result will follow.

Let's prove our claim: Let $\alpha \in A_\mu$ be a regular value
for $\Pi_\mu$. Then for any $(\alpha,u) \in M_\mu ,\
D\Pi_\mu:T_{(\alpha,u)}M_\mu \rightarrow T_\alpha A_\mu$ will be
onto. Since $\Pi_\mu$ is Fredholm of index 0, $D\Pi_\mu$ is an
isomorphism. Now consider $D\Pi_\delta:T_{(\alpha,u)}M_\delta
\rightarrow T_\alpha A_\delta$. The image of $D\Pi_\delta$ is
closed as $D\Pi_\delta$ is Fredholm. Since $T_\alpha A_\mu$ is
dense in $T_\alpha A_\delta$, and $\Pi_\mu=\Pi_\delta|_{M_\mu}$,
this would imply $D\Pi_\delta(\alpha,u)$ is also onto, and so
$\alpha$ is also regular value for $\Pi_\delta$.

Then, $M^{\alpha_0}$ is collection of isolated points and compact
in $C^{3,\delta}$ topology, and so it is finite.
\end{pf}

\section{Concluding Remarks}

\subsection{Genericity:}

These results show generic finiteness for least area planes with
$C^{3,\mu}$ smooth asymptotic boundary. But, one suspects that
this is true in general with some smoothness condition on the
boundary at infinity. In other words, it would be interesting to
show that for a class of curves in $\Si(\BH^3)$, say $C$, for any
$c\in C$, the least area planes spanning $c$ is finite.

Another interesting question in the other direction might be the
following: are there simple closed curves in $\Si(\BH^3)$, which
are spanned by infinitely many least area planes? In [A1],
Anderson constructed least area planes for any simple closed curve
in $\Si(\BH^3)$. Maybe, one can develop his methods, and find
different constructions giving different least area planes for
same curve at infinity, and then get some curve at infinity
spanned by infinitely many least area planes. In [G], Gabai gives
a relatively different construction for laminations of least area
planes for a given simple closed curve. By using that
construction, it might be possible to answer positively the
question.

\subsection{Regularity:}
In this paper, the $C^{3,\mu}$ condition can be relaxed to
$C^{2,\mu}$ except Theorem 3.5, where we use [ADN] results. It
might be interesting to know that what is the best class of curves
to achieve such a generic finiteness result. On the other hand, we
worked with Holder spaces in this paper. But, these methods can
also work for Sobolev spaces instead of Holder spaces.

\subsection{Minimal vs. Least Area:}
We could also have a generic finiteness result for minimal planes
if we can show that Lemma 4.6 is true for minimal planes, too. In
that lemma, we needed least area assumption to apply [HL] and [L]
in our situation. If one can show that same results are still true
for minimal planes (mean curvature 0 planes), then this directly
implies generic finiteness of minimal planes, by using our
methods.

\end{document}